\documentclass[letterpaper,11pt]{article}
\usepackage[latin1]{inputenc} 
\usepackage[top=3cm, bottom=3cm, left=3cm, right=3cm]{geometry}
\usepackage{graphicx,rotating,amsbsy,amssymb,amsmath,amsthm, hyperref, natbib, bbm, color}

\theoremstyle{definition}

\newtheorem{example}{Example}

\begin{document}

\title{\Large A note on recent criticisms to Birnbaum's theorem}
\author{V\'ictor Pe\~na, James O. Berger}

\maketitle 

\begin{abstract}
In this note, we provide critical commentary on two articles that cast doubt on the validity and implications of Birnbaum's theorem: Evans (2013) and Mayo (2014). In our view, the proof is correct and the consequences of the theorem are alive and well.
\end{abstract}

\section{Introduction} \label{sec:introduction}

Birnbaum's theorem \citep{birnbaum1962foundations} states that two statistical principles that are intuitively reasonable, the weak conditionality principle and the sufficiency principle, imply the likelihood principle, which is violated by statistical procedures such as $p$-values or reference priors. Ever since the result was published, there has been a lively discussion on its validity and implications. The monograph \cite{berger1988likelihood} contains a defense of the likelihood principle and responses to criticisms up to the date it was published, but the flow of articles has not stopped in the fields of statistics and philosophy of science (for example, \cite{helland1995simple}, \cite{bjornstad1996generalization}, \cite{robins2000conditioning}, \cite{sweeting2001coverage}, \cite{wechsler2008birnbaum},  \cite{grossman2011likelihood}, \cite{gandenberger2014new}). Somewhat recently, the articles \cite{evans2013does} and \cite{mayo2014} have cast doubt on the validity and implications of Birnbaum's theorem, and the goal of this note is to review and discuss their content.  \\



First, we introduce our basic notation and definitions for statistical experiment, inference bases, and informative inference: 
\begin{itemize}
\item \textbf{Statistical experiment:} A triplet $E = \{\mathcal{X}_E, \Theta_E, p_{\theta,E}\}$, where $\mathcal{X}_E$ is the sample space of the experiment, $\Theta_E$ is the parameter space, and $p_{\theta,E}$ is the sampling distribution of $E$ for $\theta \in \Theta$. As it is usual in the literature, we avoid measure-theoretical details by considering experiments with a discrete support (see Section 3.4 in \cite{berger1988likelihood} for generalizations).
\item \textbf{Inference base:} A tuple $(E,x)$ where $E$ is a statistical experiment and $x \in \mathcal{X}_E$ is an outcome from $E$.
\item \textbf{Informative inference:} $\mathbf{Ev}(E,x)$ is the informative inference (or conclusion) made by an agent given $(E,x)$. If $\mathcal{I}$ is the space of inference bases, one can think of Ev as a function from $\mathcal{I}$ to a set $\mathcal{D}$ of possible inferences. 
\item \textbf{Inferentially equivalent:} Two inference bases $(E,x)$ and $(E',x')$ are inferentially equivalent if $\mathbf{Ev}(E,x) = \mathbf{Ev}(E',x')$ (the same inferences are made given $(E,x)$ or $(E',x')$).
\end{itemize} 
Given the definitions above, we define the statistical principles at stake: the weak conditionality principle (WCP), ancillarity principle (AP), sufficiency principle (SP), and the likelihood principle (LP): 
\begin{itemize}
\item \textbf{Weak Conditionality Principle (WCP):} Consider the statistical experiments $E_1 = (\mathcal{X}_{E_1}, \Theta, p_{\theta,E_1})$, $E_2 = (\mathcal{X}_{E_2}, \Theta, p_{\theta,E_2})$ and a 50-50 mixture between $E_1$ and $E_2$, which we denote $E_{\mathrm{mix}}$. Conceptually, one can imagine that a fair coin is tossed: if it lands heads, $E_1$ is performed; if it lands tails. $E_2$ is performed. Formally, the outcome of the mixture experiment will be a pair $(j,x),$ where $j$ indicates the experiment that was performed ($j = 1$ if $E_1$ was performed, and $j=2$ if $E_2$ was performed instead), and $x \in \mathcal{X}_{E_1} \cup \mathcal{X}_{E_2}$ is the outcome of the experiment that was performed. WCP states that the informative inference given $(E_{\mathrm{mix}}, (j, x))$ from the mixture experiment should be equal to the informative inference given the inference base of the component experiment $(E_j, x)$; that is, $\mathbf{Ev}(E_{\mathrm{mix}}, (j, x)) = \mathbf{Ev}(E_j, x)$.
\item \textbf{Ancillarity Principle (AP):} Let $U$ be an ancillary statistic for $\theta$ (the distribution of $U$ does not depend on $\theta$) for which the value $u$ is observed. Then, $\mathbf{Ev}(E, (u,x)) = \mathbf{Ev}(E_{\mid U=u}, x)$, where the sampling distribution associated with $E_{\mid U=u}$ is $p_{\theta, E \mid U=u}(\cdot)$ (the conditional probability mass function of $x$ given $U=u$). In words, the ancillarity principle states that conditioning on an ancillary statistic should not change our informative inference. This principle is also known as the (strong) conditionality principle. Clearly, the selection of the component in WCP is an example of an ancillary statistic, so AP implies WCP.
\item \textbf{Sufficiency Principle (SP):} If $(E,x)$ and $(E,x')$ are such that $T(x) = T(x')$ for a sufficient statistic $T$ for $\theta$, then $\mathbf{Ev}(E, x) = \mathbf{Ev}(E,x')$.
\item \textbf{Likelihood Principle (LP):} If $(E,x)$ and $(E',x')$ are such that $p_{\theta,E}(x) = c \, p_{\theta, E'}(x')$ for $c > 0$ that does not depend on $\theta$, then $\mathbf{Ev}(E,x) = \mathbf{Ev}(E,x')$.
\end{itemize}
In our framework, $\mathbf{Ev}$ can be any function from the space of inference bases to inferences, and the mathematical role of statistical principles is restricting the set of functions that one is allowed to use. As explained in more detail in Section~\ref{sec:evans}, Evans' objections arise because a map $\mathbf{Ev}$ is not introduced. Conversely, in Section~\ref{sec:mayo} we show that the definition of the sufficiency principle in \cite{mayo2014} is different from SP (as defined in the paragraph above) and blocks Birnbaum's proof.

\section{Evans' objections}
\label{sec:evans}
Evans defines the statistical principles as the following set relations on $\mathcal{I} \times \mathcal{I}$:\footnote{We use a slightly different notation than in \cite{evans2013does}. We define $C$ as a formalization of WCP and $A$ as a formalization of AP. However, \cite{evans2013does} does not consider WCP at all and defines a set relation (which is denoted $C$ in Evans' article) which is equivalent to our $A$. We apologize for the possible confusion that this change might cause.}
\begin{itemize}
\item $C$: $(E,x) \sim_{C} (E',x')$ if and only if $E = E_{\mathrm{mix}}$, $x = (j,x_j)$, $E' = E_j$, and $x' = x_j$ as in the definition of WCP in Section~\ref{sec:introduction} (or with roles of $(E,x)$ and $(E',x')$ reversed).
\item $A$:  $(E,x) \sim_{A} (E',x')$ if and only if $x = (u,x')$ and $E' = E_{|U=u}$, where $U = u$ and $E_{|U=u}$ are as defined for AP in Section~\ref{sec:introduction} (or with roles of $(E,x)$ and $(E',x')$ reversed). 
\item $S$:  $(E,x) \sim_{S} (E',x')$ if and only if there exists a sufficient statistic $T$ for $\theta$ such that $T(x) = T(x')$.
\item $L$:  $(E,x) \sim_{L} (E',x')$ if and only if $p_{\theta, E}(x) = c \, p_{\theta, E'} (x')$ for a constant $c > 0$ which does not depend on $\theta$.
\end{itemize}
This approach is different from the one taken in Section~\ref{sec:introduction} and the one in \cite{birnbaum1962foundations} because $\mathbf{Ev}$ is not defined or used at all in the definitions. Nonetheless, the set relations are very similar to the principles defined in Section~\ref{sec:introduction}: they are of the form $(E,x) \sim_{P} (E',x')$ if and only if $\mathbf{Ev}(E,x) = \mathbf{Ev}(E',x')$ by an application of a principle statistical $P$. According to Evans, ``A basic step missing in Birnbaum (1962) was to formulate the principles as relations on the set $\mathcal{I}$ of all model and data combinations.'' But the definition of a function $\mathbf{Ev}$ automatically induces an equivalence relation on $\mathcal{I} \times \mathcal{I}$ (the kernel of $\mathbf{Ev}$): $(E,x) \sim (E,x')$ if and only if $\mathbf{Ev}(E,x) = \mathbf{Ev}(E',x')$. If we accept WCP and SP as defined in Section~\ref{sec:introduction}, the equivalence relation on $\mathcal{I} \times \mathcal{I}$ induced by accepting WCP and SP implies that bases with proportional likelihoods are equivalent because they map to the same value. \\

Evans shows that statistical principles defined as set relations need not be equivalence relations (for instance, $A$ and $C$ as defined above are not). Even if two statistical principles are equivalence relations, their union may not be because it could fail to be transitive: if $P_1$ and $P_2$ are set relations formalizing statistical principles, it is possible that the inference bases $(E_1, x_1)$ and $(E_2, x_2)$ are inferentially equivalent with respect to $P_1$ and $(E_2, x_2)$ and $(E_3, x_3)$ are inferentially equivalent according to $P_2$ but $(E_1, x_1)$ and $(E_3, x_3)$ are not inferentially equivalent according to either $P_1$ or $P_2$ alone. \\ 

Birnbaum's argument is a neat (and in our view, transparent) illustration of this phenomenon: the inference bases with proportional likelihoods are shown to be equivalent by a chain of applications of WCP and SP, but they are not equivalent according to either WCP or SP individually. This implies that $L \neq S \cup C$. The correct result is that $L$ is equal to the smallest equivalence relation generated by $S \cup C$, and Evans argues that extending statistical principles that are originally defined as set relations to equivalence relations requires further justification. \\

Here is a simple one: if we define the sufficiency principle and the weak conditionality principle as the set relations $S$ and $C$ and introduce a function $\mathbf{Ev}$ with the minimal requirement that $\mathbf{Ev}(E,x) = \mathbf{Ev}(E',x')$ if and only if $(E,x) \sim_{C} (E',x')$ or $(E,x) \sim_{S} (E',x')$, the equivalence relation on $\mathcal{I} \times \mathcal{I}$ generated by $\mathbf{Ev}$ is precisely the smallest equivalence relation generated by $S \cup C$, which in this case is $L$. In general, if we define statistical principles $P_1, P_2, \, ... \, , P_k$ as set relations on $\mathcal{I} \times \mathcal{I}$ and introduce $\mathbf{Ev}$ with the property $\mathbf{Ev}(E,x) = \mathbf{Ev}(E',x')$ if and only if $(E,x) \sim_{P_i} (E',x')$ for some $i \in \{1,2, \, ... \, k\}$, the equivalence relation on $\mathcal{I} \times \mathcal{I}$ induced by $\mathbf{Ev}$ is equal to the smallest equivalence relation generated by $P_1, P_2, \, ... \, , P_k$. Defining statistical principles as set relations on $\mathcal{I} \times \mathcal{I}$ and introducing $\mathbf{Ev}$ as we just did is equivalent to stating the definitions in terms of $\mathbf{Ev}$ in the first place as in Section~\ref{sec:introduction}. \\

Since the notation $\mathbf{Ev}$ is very explicit in \cite{birnbaum1962foundations}, we believe that the definition of the principles in terms of set relations was implied by the fact that $\mathbf{Ev}$ is a function. But even within a framework where $\mathbf{Ev}$ is not defined, the smallest equivalence relation generated by a collection of principles has a straightforward interpretation: its elements are (exclusively) the result of a chain of applications of the principles we wish to respect. Rejecting the extension implies rejecting the equivalence of inference bases that can be shown to be equivalent by a number of applications of our principles. We believe, then, that the extension is also justified if $\mathbf{Ev}$ is not introduced.  \\

Now we turn to an example in \cite{evans2013does} that shows that $A$ is not transitive and illustrates some of the issues that were commented in the paragraphs above.

\begin{example} \label{ex:evans} (\cite{evans2013does}, pg. 2651)
Let $\mathcal{X}_E = \{1,2\}\times\{1,2\}$, $\Theta_E = \{1,2\}$, with $p_{E,\theta}$ given in Table~\ref{tab:unconditional}. Both $U(x_1, x_2)= x_1$ and $V(x_1,x_2)=x_2$ are ancillary, and the conditional models upon observing $U = 1$ and $V = 1$ are given in Tables~\ref{tab:condu} and~\ref{tab:condv}. This example shows that $A$ is not transitive: $(E,(x_1,x_2)) \sim_{A} (E_{\mid U}, x_2)$ and $(E,(x_1,x_2)) \sim_{A} (E_{\mid V}, x_1)$, but $(E_{\mid U}, x_2) \not \sim_A (E_{\mid V}, x_1)$ because there is no ancillary statistic linking the two conditional models. However, using the definitions in Section~\ref{sec:introduction} (or equivalently, using $A$ and introducing Ev with the property $\mathbf{Ev}(E,x) = \mathbf{Ev}(E',x')$ if and only if $(E,x) \sim_{A} (E',x')$), we have $\mathbf{Ev}(E,(x_1,x_2)) = \mathbf{Ev}(E_{\mid U}, x_2)$ and $\mathbf{Ev}(E,(x_1,x_2)) = \mathbf{Ev}(E_{\mid V}, x_1)$, so $ \mathbf{Ev}(E_{\mid U}, x_2) =  \mathbf{Ev}(E_{\mid V}, x_1)$.
\begin{table}[h!]
\centering
    \caption{Unconditional model (rows: sampling distributions for $\theta \in \{1,2\}$)}
    \label{tab:unconditional}
    \begin{tabular}{l|llll}
    $(x_1,x_2)$         & $(1,1)$ & $(1,2)$ & $(2,1)$ & $(2,2)$ \\
    \hline
    $f_{E,\theta=1}(x_1, x_2)$ & 1/6     & 1/6     & 2/6     & 2/6     \\
    $f_{E,\theta=2}(x_1,x_2)$  & 1/12    & 3/12    & 5/12    & 3/12    \\
    \end{tabular}
\end{table}

\begin{table}[h!]
\centering
    \caption{Conditional model when $U =1$  (rows: sampling distributions for $\theta \in \{1,2\}$)}
    \label{tab:condu}
    \begin{tabular}{l|llll}
    $(x_1,x_2)$               & $(1,1)$ & $(1,2)$ & $(2,1)$ & $(2,2)$ \\
    \hline
    $f_{E,\theta=1}(x_1, x_2 \mid U=1)$ & 1/2     & 1/2     & 0       & 0       \\
    $f_{E,\theta=2}(x_1,x_2 \mid U=1)$  & 1/4     & 3/4     & 0       & 0       \\
    \end{tabular}
\end{table}

\begin{table}[h!]
\centering
    \caption{Conditional model when $V =1$  (rows: sampling distributions for $\theta \in \{1,2\}$)}
    \label{tab:condv}
    \begin{tabular}{l|llll}
    $(x_1,x_2)$               & $(1,1)$ & $(1,2)$ & $(2,1)$ & $(2,2)$ \\
    \hline
    $f_{E,\theta=1}(x_1, x_2 \mid V=1)$ & 1/3     & 0       & 2/3     & 0       \\
    $f_{E,\theta=2}(x_1,x_2 \mid V=1)$  & 1/6     & 0       & 5/6     & 0       \\
    \end{tabular}
\end{table}
\end{example}
Quoting \cite{evans2013does}: ``Saying that such models [the conditional models in Tables~\ref{tab:condu},~\ref{tab:condv}] contain an equivalent amount of statistical information is clearly a substantial generalization of [$A$]. To measure the accuracy of this estimate we can compute the conditional probabilities based on the two inference bases, namely,
$$ 
\mathbb{P}_{\theta=1}(\widehat{\theta} = 1 \mid U = 1) = 1/2, \, \, \, \mathbb{P}_{\theta=2}(\widehat{\theta} = 1 \mid V = 1) = 3/4
$$
and so the accuracy of $\widehat{\theta}$ is quite different depending on whether we [condition on $U$ or $V$]. It seems unlikely that we would interpret these inference bases as containing an equivalent amount of information in a frequentist formulation of statistics.'' \\

Concluding that the inference bases are equivalent with respect to $A$ is a consequence of introducing $\mathbf{Ev}$ with the property $\mathbf{Ev}(E,x) = \mathbf{Ev}(E',x')$ if and only if $(E,x) \sim_A (E',x')$. Also, the likelihood ratio of $\theta = 1$ to $\theta = 2$ equals 2 if we condition on either $U$ or $V$, which is unsurprising because, as Evans proves, AP equals LP. We agree with Evans in that this example shows that accepting AP can be problematic for frequentist statisticians: there are two ancillary statistics we can condition on, there is no apparent reason one should prefer one over the other and, unfortunately, standard errors and $p$-values depend depend on the choice of ancillary. We return to this point in Section~\ref{sec:ancillaries}.  \\

After showing that AP is equivalent to LP, Evans concludes that SP is redundant in Birnbaum's argument. Then, Example~\ref{ex:evans} leads him to cast doubt on the impact of Birnbaum's result because he believes that many statisticians would not accept AP (or equivalently, $A$ and the equivalences generated by the principle). But SP is certainly not redundant if only WCP is assumed (recall that WCP only requires equivalence of 50-50 mixtures), and WCP and SP also imply LP. In Example~\ref{ex:evans}, the conditional experiments are not equivalent according to WCP, and the smallest equivalence relation containing $C$ would only add cases where mixture experiments with different components (or different probabilities of performing them) were considered, but the same component experiment was performed and the same result was obtained. Finally, we agree with Evans that accepting statistical principles may induce unexpected equivalences between inference bases, which is precisely what makes Birnbaum's result surprising and relevant.
 
\section{Mayo's objections} \label{sec:mayo}
In our view, the objections to Birnbaum's proof in \cite{mayo2014} stem from using a definition for the sufficiency principle that is different from that in Section~\ref{sec:introduction}. 
We believe that introducing new notation that makes an explicit distinction between the output of methods and the inference made by an agent that is using them is helpful for understanding the arguments:
\begin{itemize}
\item $\mathbf{M}(E, x)$: Result of a applying a method $M$ to the inference base $(E,x)$.
\item $\mathbf{Ev}(E,x)$: Inference made by an agent given $(E,x)$ (as in Section~\ref{sec:introduction}).
\end{itemize}
Given $(E,x)$, the agent makes informative inferences $\mathbf{Ev}(E.x)$ by means of $\mathbf{M}(E',x')$ for some $(E',x')$ which may not be equal to $(E,x)$. The interpretation of $\mathbf{M}(E,x) = \mathbf{M}(E',x')$ is that the ``output'' of applying a method $\mathbf{M}$ to $(E,x)$ and $(E',x')$ is the same (one can imagine that $\mathbf{M}$ is a function in some programming language that takes $E$ and $x$ as inputs), whereas $\mathbf{Ev}(E,x) = \mathbf{Ev}(E',x')$ means that an agent makes the same informative inferences given $(E,x)$ or $(E',x')$. This distinction is somewhat obscured in \cite{mayo2014}, as she defines
\begin{itemize}
\item $\mathrm{Infr}_E[x]$: The parametric statistical inference from a given or known $(E,z)$.
\item $(E',x') \Rightarrow \mathrm{Infr}_E[x]$: An informative parametric inference about $\theta$ from given $(E,x)$ is to be computed by means of $\mathrm{Infr}_E[x]$.
\end{itemize}
The definition of $\mathrm{Infr}_E[x]$ and the name ``Infr'' suggest that $\mathrm{Infr}_E[x] = \mathbf{Ev}(E,x)$. However, the second definition implies that $\mathrm{Infr}_E[z]$ need not be equal to the final inference $\mathbf{Ev}(E,x)$. This is explicit in her definition of the weak conditionality principle (WCP):
\begin{itemize}
\item \textbf{WCP:} Given $(E_\mathrm{mix}, (j,x_j))$, condition on the $E_j$ producing the result: $(E_\mathrm{mix}, (j,x_j)) \Rightarrow \mathrm{Infr}_{E_j}[x_j]$. Do not use the unconditional formulation: $(E_\mathrm{mix}, (j,x_j)) \not \Rightarrow \mathrm{Infr}_{E_\mathrm{mix}}[(j,x_j)]$.
\end{itemize}
Using our notation, this definition is equivalent to the WCP in Section~\ref{sec:introduction}. However, Mayo defines the sufficiency principle as follows
\begin{itemize}
\item \textbf{SP2}: If there exists a sufficient statistic $T$ for $\theta$ and $T(x) = T(x')$, then $\mathrm{Infr}_E[x] = \mathrm{Infr}_E[x']$,
\end{itemize}
which is different from SP, and can be recast as 
\begin{itemize}
\item \textbf{SP2:} If there exists a sufficient statistic $T$ for $\theta$ and $T(x) = T(x')$, then $\mathbf{M}(E, x) = \mathbf{M}(E, x')$.
\end{itemize}
The key point is that WCP is a property of $\mathbf{Ev}$ and SP2 is a property of $\mathbf{M}$. If this distinction is made, LP does not follow. The distinction between $\mathbf{Ev}$ and $\mathbf{M}$ is not made in \cite{birnbaum1962foundations} or Section~\ref{sec:introduction}. The following example, which is a slight modification of the example presented in Section 4. in \cite{mayo2010iii}, puts the notation in context and makes clear why WCP and SP2 do not imply LP.

\begin{example} \label{ex:mayo}
Consider binomial and negative binomial experiments $$E_1 = \{ \{0,1,2, \, ... \, , n\}, \, \Theta, \, \mathrm{Binomial}(n, \theta) \}, \qquad E_2 = \{ \{0,1,2, \, .... \}, \, \Theta, \, \mathrm{NegBinomial}(k, \theta) \}.$$ Suppose that a fair coin is flipped and $E_1$ is performed if the coin lands heads and $E_2$ is performed if it lands tails. Let $E_{\text{mix}}$ denote the ``mixture'' experiment. The outcome of $E_{\text{mix}}$ is $(j, x)$, with $j \in \{1,2\}$ ($j=1$ if $E_1$ is performed and $j=2$ if $E_2$ is performed) and $x = (k, n-k)$, where $k$ and $n-k$ are the number of successes and failures observed after performing $E_j$. The statistical method ${M}(E,x)$ is the one-sided $p$-value for testing $\theta = \theta_0$ against $\theta > \theta_0$:
\begin{align*}
 \mathbf{M}(E_1, x) &= \mathbb{P}( \mathrm{Binomial}(n,\theta_0) \ge x) \\
 \mathbf{M}(E_2, x) &= \mathbb{P}( \mathrm{NegBinomial}(r,\theta_0) \ge x) \\
 \mathbf{M}(E_{\mathrm{mix}}, x) &= 0.5 \,  \mathbb{P}( \mathrm{Binomial}(n,\theta_0) \ge x) + 0.5 \, \mathbb{P}( \mathrm{NegBinomial}(k,\theta_0) \ge x).
\end{align*}
We assume that the agent makes inference using the rule $\mathbf{Ev}(E_{\mathrm{mix}}, (j,x)) = \mathbf{M}(E_{j}, x)$. The statistic $T(j,x) = (1,x)$ is sufficient for $\theta$ with respect to $E_{\text{mix}}$, and it satisfies both $T(1,x) = T(2,x)$ and $\mathbf{M}(E_{\text{mix}},(1,x)) = \mathbf{M}(E_{\text{mix}},(2,x))$, so SP2 is respected. WCP is automatically satisfied because the inference rule is $\mathbf{Ev}(E_{\text{mix}}, (j,x)) =  \mathbf{M}(E_{j}, x) = \mathbf{Ev}(E_j, x_j)$ (the inference rule is chosen so that WCP is respected). It follows that WCP and SP2 do not imply LP. 
\end{example}
According to the definitions in \cite{mayo2014}, WCP and SP2 do not imply LP, as seen in the example above. Where does Birnbaum's proof go wrong? With WCP as stated, the mixture experiments are inferentially equivalent to the performed components: $\mathbf{Ev}(E_{\mathrm{mix}}, (1,x_1)) = \mathbf{Ev}(E_1, x_1)$ and $\mathbf{Ev}(E_{\mathrm{mix}}, (2,x_2))$ $= \mathbf{Ev}(E_2, x_2)$. However, SP2 does not imply $\mathbf{Ev}(E_{\mathrm{mix}}, (1,x_1))$ $= \mathbf{Ev}(E_{\mathrm{mix}}, (2,x_2))$: instead, it requires $\mathbf{M}(E_{\mathrm{mix}}, (1,x_1)) = \mathbf{M}(E_{\mathrm{mix}}, (2,x_2))$. However, $\mathbf{Ev}(E_{\mathrm{mix}}, (1,x_1))$ need not be equal to $\mathbf{Ev}(E_{\mathrm{mix}}, (2,x_2))$. These definitions allow Mayo to claim that, in Example~\ref{ex:mayo}, reporting the conditional $p$-value according to the sampling distribution of the component experiment that was performed does not violate the sufficiency principle. In contrast, reporting the conditional $p$-value is a violation of SP as defined in Section~\ref{sec:introduction} (and the proof of WCP and SP implies LP goes through as usual). Critically, note that SP states that if there exists a sufficient statistic, the inferences bases are inferentially equivalent, but there is no requirement that said sufficient statistic be used for our final inferences. If that were the case, it would imply that SP instructs to use of the unconditional $p$-value. The reason that Birnbaum's proof does not go through in this framework hinges on the distinction of $\mathbf{Ev}$ and $\mathbf{M}$: if we define a new WCP2 as a property of $\mathbf{M}$ (so that both SP2 and WCP2 were properties of $\mathbf{M}$), reporting the conditional $p$-value in Example~\ref{ex:mayo} would violate WCP2, as $\mathbf{M}(E_j, x_j) \neq \mathbf{M}(E_{\mathrm{mix}},(j,x_j))$ (and WCP2 and SP2 would, of course, imply a version of LP written in terms of $\mathbf{M}$). 

\section{Can AP be applied in frequentist statistics?} \label{sec:ancillaries}

We briefly discuss the applicability of the ancillarity principle in frequentist inference, motivated by comments in \cite{coxmayo2010}, \cite{evans2013does}, and \cite{mayo2014}. Since AP is equivalent to LP, frequentist statisticians that want to make conditional frequentist statements have to propose restricted versions of AP. Additionally, it is of utmost importance to find well-defined criteria for choosing among ancillaries because, as we have seen in Example~\ref{ex:evans}, there are instances where there are multiple ancillaries one can condition on that give rise to different conditional $p$-values or standard errors. Some authors have proposed restricting the set of ancillaries to condition on (\cite{durbin1970birnbaum}, \cite{kalbfleisch1975}), but this approach is problematic because there are examples where several ancillaries satisfy the restrictions (see \cite{basu1964recovery} for examples and \cite{dawid2011basu} for a concise and lucid review on the ancillarity principle and the issues that have been mentioned in this paragraph). To the best of our knowledge, there is no (restricted) formulation of AP that instructs which ancillary one should use for any given problem, and as a result is no adequate definition for a restricted ancillarity principle that is not equivalent to the likelihood principle (\cite{cox1971choice} provides a heuristic that works when applied to an example in \cite{basu1964recovery}, but it does not give a definite answer in other problems and it is not regarded as a general solution to this problem). Another issue is that there are examples where a conditional analysis is clearly desirable, but useful ancillary statistics are not available. We present two examples below. 

\begin{example} \label{ex:example8} (Example 8 in \cite{berger1988likelihood}) 
Let $\Theta = [0,1)$, $P(X = \theta) = 1-\theta$, and $P(X = 0) = \theta$. Consider the confidence set $C = \{X\}$. Unconditionally, $P(\theta \in C) = 1-\theta$. However, if $X > 0$, we know that $C = \{X\}$ contains $\theta$ with probability 1, but $X$ is not ancillary, so the ancillarity principle would not allow conditioning on its value.
\end{example}

\begin{example} \label{ex:noancillary}
Let $X_1,X_2$ be independent and identically distributed random variables with $P(X_i = \theta-1) = P(X_i = \theta+1) = 1/2$ for $i \in \{1,2\}$. Let $D = |X_1-X_2|/2$, which is ancillary with $P(D=1) = P(D=0) = 1/2$. Suppose we want to evaluate the quality of the estimator $T = X_{(1)}+1$ ($X_{(1)}$ is the minimum of $X_1$ and $X_2$). Conditioning on $D$, we know that $P(T = \theta \mid D = 1) = 1$ and $P(T = \theta \mid D = 0) = 1/2$, and \cite{coxmayo2010} would propose reporting inferences conditional on $D$ because it is more informative than an unconditional analysis. But now consider the following modification: $P(X_i = \theta+1) = 1/2+\theta \epsilon$ and $P(X_i = \theta -1) = 1/2-\theta \epsilon$ for a known $\epsilon \in [0, 1]$ and $\theta \in [-1/(2\epsilon),1/(2\epsilon)]$. The original example is a particular case with $\epsilon = 0$. If $\epsilon \neq 0$, $D$ is not ancillary anymore, despite the fact that if $\epsilon$ is small (say $\epsilon = 10^{-100}$) we are essentially in the same situation as if $\epsilon = 0$. In addition, if $\epsilon \neq 0$, there are (even more) cases where we can retrieve $\theta$ with probability 1 given the data. Indeed, if $X_1 \neq X_2$ we still have that $\theta = X_{(1)}+1$, but now there are cases where we know the value of $\theta$ exactly even if $X_{(1)} = X_{(2)}$. Let $A_{\theta-1} = [-1/(2\epsilon)-1, 1/(2\epsilon)-1]$ and $A_{\theta+1} = [-1/(2\epsilon)+1,1/(2\epsilon)+1]$. If $X_{(1)} \in A_{\theta-1} \setminus A_{\theta+1}$, then $\theta = X_{(1)}+1$; analogously, $\theta = X_{(1)}-1$ whenever $X_{(1)} \in A_{\theta+1} \setminus A_{\theta-1}$. Note that if $\epsilon > 1/2$, $A_{\theta-1} \cap A_{\theta+1} = \emptyset$ and we can always retrieve the value of $\theta$. If we want to assess the performance of $T$ conditionally, we know that
\begin{align*}
P(T &= \theta \mid X_{(1)} \neq X_{(2)}) = 1 \\
P(T  &= \theta \mid X_{(1)} = X_{(2)}, X_{(1)} \in (A_{\theta-1}\setminus A_{\theta+1})) = 1 \\
P(T &= \theta \mid X_{(1)} = X_{(2)}, X_{(1)} \in (A_{\theta+1}\setminus A_{\theta-1})) = 0 \\
P(T &= \theta \mid X_{(1)} = X_{(2)}, X_{(1)} \in A_{\theta-1} \cap A_{\theta+1}  )= \frac{(1/2-\theta \epsilon)^2}{(1/2-\theta \epsilon)^2+(1/2+\theta \epsilon)^2},
\end{align*}
but unconditionally $P(T = \theta) = 1 - (1/2+\theta \epsilon)^2$, which depends on $\theta$ and ranges from 1 to 0 for $\theta \in [-1/(2\epsilon),1/(2\epsilon)]$. Therefore, the confidence level of the set $C = \{T\}$ is $\inf P_\theta(\theta \in C) = 0$, which is clearly undesirable and misleading (especially in cases where $\epsilon > 1/2$, where a conditional analysis reveals if $T = \theta$ with probability 0 or 1 depending on the data). As an aside, a modified estimator that takes on the value $X_{(1)}-1$ whenever $X_{(1)} = X_{(2)}$ and $X_{(1)} > 0$ has better performance, but we used $T$ for illustrative purposes.
\end{example}

Finally, we note that applying the ancillarity principle can be suboptimal according to strictly frequentist criteria: in practice, there are cases where an unconditional test is preferable to a conditional test, as in the following example inspired by an Example in \cite{cox1958some}.

\begin{example}
Suppose a production line is periodically tested to see if it is operating correctly. If correct, it produces a part of diameter 1. Periodically it goes out of line and then produces parts with diameter 1.1. In the testing, the parts are measured with one of two measuring instruments, an old one which produces a normal observation with mean the true diameter of the part and standard deviation 0.1, and a new measuring instrument which produces a normal observation with mean the true diameter and standard deviation 0.05. The old and new measuring instruments are each available with probability 1/2 (as there is another production line for which they are also used). If the production line is deemed to be out of line, it must be shut down and reset, at considerable expense. The company does a cost-benefit analysis and determines that it will be optimal to control overall Type I error in the testing at the 0.05 level. This is a scenario in which frequentist analysis is absolutely appropriate, in that there is true long-term repetition of the test. Also, the cost-benefit analysis is presumably carried out in a Bayes-frequentist sense, since historical levels of in-line and out-of-line must be taken into account. If the company followed WCP, they would do the 0.05 level test conditional on which measuring instrument is being used at each test. But this will lose the company money, as the power of this test for detecting an out-of-line process (which is 0.646) is 9\% less than that of the most powerful test (which is 0.694). This most powerful test corresponds to using Type I error probabilities of 0.099 and 0.001 for the old and new measuring instruments, respectively.
\end{example}

The example above is interesting in that it suggests that, for frequentists, the only way to implement the conditionality principle is to use a method that is compatible with Bayesian reasoning (as the unconditional test would be equivalent to the Bayes rule with respect to the loss function implied by the cost-benefit analysis). This is not surprising, given the complete class theorems that show that optimal frequentist decision procedures are necessarily Bayesian. 

\section{Conclusions}

The articles \cite{evans2013does} and \cite{mayo2014} contain thought-provoking discussions about the conditions under which the result in \cite{birnbaum1962foundations} is valid, but that neither of them show that WCP and SP do not imply LP according to the definitions in Section~\ref{sec:introduction}, which, in our view, are equivalent to the definitions in \cite{birnbaum1962foundations}.  \\

\newpage

Evans avoids introducing Ev, which is central in Birnbaum's argument, and defines statistical principles as set relations on the (product) space of inferences. If $\mathbf{Ev}$ is introduced with the property $\mathbf{Ev}(E,x)$ = $\mathbf{Ev}(E',x')$ if and only if $(E,x) \sim_C (E',x')$ or $(E,x) \sim_S (E',x')$  Birnbaum's result follows. If we stick to Evans' framework, the union of the set relation defined by the sufficiency principle ($S$) and the conditionality principle ($C$) does not equal the set relation defined by the likelihood principle ($L$). This result might seem surprising at first glance but, if it were true, two inference bases with proportional likelihoods would be equivalent according to either the sufficiency principle or the conditionality principle individually, which is clearly false. What is true is that the smallest equivalence relation generated by $S \cup C$ equals $L$. As explained in Section~\ref{ex:evans}, the equivalence relation generated by $S \cup C$ only contains inference bases that are equivalent to a chain of applications of the principles.  \\

Mayo defines statistical principles making a distinction between the output of methods ($\mathbf{M}$) and the inferences that are made by an agent using them ($\mathbf{Ev}$): the weak conditionality principle is defined as a property of $\mathbf{Ev}$, whereas the sufficiency principle is defined as a property of $\mathbf{M}$. For example, this distinction allows Mayo to claim that in a mixture experiment where a Negative Binomial or Binomial experiment is selected with equal probability, reporting the conditional $p$-value does not result in a violation of the sufficiency principle (see Example~\ref{ex:mayo}). In the framework of \cite{mayo2014}, the weak conditionality principle and the sufficiency principle do not imply the likelihood principle, but the definition of the sufficiency principle differs from that in \cite{birnbaum1962foundations} because the distinction between the ouput of methods and informative inference is not made.

\bibliographystyle{chicago}  
\bibliography{objectionslp} 

\end{document}